\documentclass[compress, number]{elsarticle}

\usepackage{hyperref}
\hypersetup{
    colorlinks=true,
    linkcolor=blue,
    urlcolor=cyan,
    citecolor=blue
    }
\usepackage{bm}
\usepackage{amsmath}
\usepackage{amssymb}
\usepackage{enumitem}

\newcommand{\dif}{\mathcal{M}}
\newcommand{\dar}{\mathcal{D}}
\newcommand{\ten}{\mathbf}
\newcommand{\gt}{\bm}
\newcommand{\scalar}[1]{\langle#1\rangle}

\newcommand{\norm}[1]{\lVert#1\rVert}
\newcommand{\str}{\varepsilon}
\newcommand{\dd}{\mathrm{d}}

\renewcommand{\[}{\begin{equation}}
\renewcommand{\]}{\end{equation}}
\newcommand{\R}{\mathbb{R}}
\let\div\relax
\DeclareMathOperator{\div}{div}
\DeclareMathOperator{\adj}{adj}

\newtheorem{example}{\bf Example}
\newtheorem{lemma}{\bf Lemma}
\newtheorem{definition}{\bf Definition}
\newtheorem{theorem}{\bf Theorem}
\newtheorem{corollary}{\bf Corollary}
\newproof{proof}{Proof}

\begin{document}
\begin{frontmatter}
    
\title{Effective isometries of periodic shells}

\author{Hussein Nassar\corref{cor1}}
\ead{nassarh@missouri.edu}
\author{Andrew Weber}

\address{Department of Mechanical and Aerospace Engineering, University of Missouri, Columbia, 65211 MO, USA}

\cortext[cor1]{Corresponding author}

\begin{abstract}
We argue that the standard classification of isometric deformations into infinitesimal v.s. finite is inadequate for the study of compliant shell mechanisms. Indeed, many compliant shells, particularly ones that are periodically corrugated, exhibit low-energy deformations that are far too large to be infinitesimally isometric and far too rich to be finitely isometric. Here, rather than abandon the geometric standpoint in favor of a full theory of elastic shells, we introduce the concept of \emph{effective isometric deformations} defined as deformations that are first-order isometric in a small scale separation parameter given by the ratio of the size of the corrugation to the size of the shell. The main result then states that effective isometries are solutions to a quasilinear second-order PDE whose type is function of an effective geometric Poisson's ratio. The result is based on a self-adjointness property of the differential operator of infinitesimal isometries; it holds for periodic surfaces that are smooth, piecewise smooth or polyhedral, i.e., periodic surfaces with or without straight and curved creases. In particular, it unifies and generalizes a series of previous results regarding the effective Poisson's ratio of parallelogram-based origami tessellations. Numerical simulations illustrate and validate the conclusions.
\end{abstract}

\end{frontmatter}


\section{Introduction}\label{sec1}
The flexure and membrane strain energies of a thin elastic shell occur in a ratio of order $(t/\rho)^2$ where $t$ is thickness and $\rho$ is a typical radius of curvature~\cite{Landau1986}. Thus, in the limit $t/\rho\to 0$, the deformation of a shell is, ideally, inextensional~\cite{Rayleigh1894}. In mathematical jargon, an inextensional deformation defines an \emph{isometry} of the shell's midsurface, namely a deformation that preserves lengths as measured on the midsurface. In the case of small deflections, the relevant notion is that of an \emph{infinitesimal isometry}, i.e., a deformation that preserves lengths up to first order in the deflection's amplitude; see~\cite{Connelly1993, Spivak1999a, Ivanova-Karatopraklieva1994, Ivanova-Karatopraklieva1995} for a comprehensive introduction to the topic. The assumption of inextesional deformations is at the basis of the \emph{flexure theory} of shells which proceeds in two steps: first, determine the isometric deformations that are compatible with the imposed boundary conditions; second, determine, among isometric deformations, the one that minimizes strain energy. Evidently, for the flexure theory to be fruitful, non-trivial isometric deformations must exist. Should they not, meaning that all deformations involve some extension, then the shell is better described by a \emph{membrane theory} where flexure energy is neglected, it being infinitesimally smaller than membrane energy. Nonetheless, this dichotomy is an oversimplification of the behavior of elastic shells, particularly so in cases where a small extension coexists with a large flexion. If so, then a complete theory of thin elastic shells, one where extension competes with bending, becomes necessary, e.g., Föppl-von Kármán theory of plates or Koiter's theory of shells~\cite{ciarlet2006, Audoly2010}.

Modern investigations carried into the design and modeling of compliant shell mechanisms provide further instances where the ``flexure v.s. membrane'' dichotomy fails~\cite{Norman2008,Norman2009,norman2009phd,Seffen2012}. Indeed, it has been observed that compliant shell mechanisms exhibit low-energy deformation modes that involve large flexure-dominant yet non-strictly inextensional deflections; see, e.g.,~\cite{Schief2008, Demaine2011}. These deformations are neither infinitesimally nor finitely isometric meaning that the use of flexure theory is not justified. At the same time, these deformations remain dominated by flexure in a way that invalidates the use of membrane theory. Short of using a complete theory of thin elastic shells where all deformations are equally likely \emph{a priori}, it is desirable to have an alternative, more insightful, description of the low-energy deformation modes, one that establishes how they are related to isometric deformations and clarifies if and how they emerge under given boundary conditions. Thus, the main purpose of the present contribution is to introduce a notion of ``approximately isometric but non-infinitesimal'' deformations specifically in the context of periodic shells. 

Previous attempts to characterize the low-energy deformation modes of compliant shell mechanisms have focused on origami tessellations. Typically, the observed deformations are explained by the emergence of ad-hoc mid-facet ``faint'' creases that enrich the space of admissible deformations to restore strict inextensibility at the cost of smoothness, see, e.g.,~\cite{Schenk2011a}. On one hand, origami models folded out of paper and thin plastic shims appear to contradict this hypothesis. On the other hand, it is unclear how to adapt the idea of emergent faint creases to shells with no planar facets, e.g., to origami tessellations with curved creases or even to smoothly corrugated shells. Here, we propose a geometric theory of low-energy deformation modes that makes no \emph{a priori} ad-hoc assumptions. The theory is applicable to smooth and piecewise smooth periodic shells, including shells with or without straight and curved creases. The key idea, inspired by homogenization theory, is to define low-energy deformations, henceforth referred to as \emph{effectively isometric deformations}, as deformations that are isometric up to first order in the periodicity wavelength, rather than in the deflection's amplitude. As mentioned, the theory remains purely geometrical; the fact that effective isometries are indeed of low energy is formally justified by the smallness of any potential extensional deformations.

The main result of the theory is a quasi-linear Partial Differential Equation (PDE) satisfied by the parametrization of the effective surface of the shell, i.e., the surface obtained by averaging out any periodic corrugations or creases. An appealing feature of origami tessellations, and compliant shell mechanisms more generally, is how approximately inextensional deformations on the level of the shell can lead to large effective extensions and curvatures on the level of the effective surface~\cite{Lebee2015}. The PDE informs precisely on that: it constrains the effective curvatures relative to the effective extensions. In particular, the type of the PDE is function of the sign of an effective Poisson's ratio defined as the ratio of principal effective strains. Moreover, the PDE states that the effective Poisson's ratio is equal to the ratio of normal curvatures in the principal directions of effective strain. This property has been demonstrated, numerically or analytically, for a number of origami and origami-like tessellations including the Miura-ori~\cite{Schenk2013, Wei2013}, the ``eggbox'' pattern~\cite{Nassar2017a}, the ``morph'' pattern~\cite{Pratapa2019}, ``zigzag sums''~\cite{Nassar2022}, and, more recently, a class of curved-crease origami tessellations~\cite{Karami2023}. The present theory unifies these instances and further shows that they are but particular cases of a more general theory valid for all periodic shells.

The first part of the paper introduces necessary preliminaries and motivates a restricted theory through a heuristic analysis of smooth surfaces that are Cartesian graphs of some function. Lemma~\ref{lem:self1} is the crucial ingredient of this restricted theory: it proves a self-adjointness property of the differential operator of infinitesimal isometries, i.e., the linearized two-dimensional Monge-Ampère operator. This property is somewhat elementary but does not appear to be well-known; thus, we take the time to introduce it and justify it. Theorem~\ref{prop:main1} is the main result of the restricted theory: it places a linear constraint on the effective curvatures and twist of smooth periodic graphs. The first part is concluded by a brief case study of graphs \emph{of translation}: Corollaries~\ref{cor:trans1},~\ref{cor:trans2} and~\ref{cor:trans3} describe how graphs of translation can or cannot bend; Example~\ref{exa:ex} describes Finite Element simulations that validate the theoretical prediction regarding the effective Poisson's ratio and the effective normal curvatures. This first part of the paper is self-sufficient and should be enough to communicate the main findings of the theory as well as its main use cases. The second part of the paper is somewhat more involved and is dedicated to the complications that arise in dealing with immersed surfaces that are not graphs and/or not smooth. It necessitates the use of two-scale asymptotic expansions; here too, a crucial ingredient is a self-adjointness property of a relevant differential operator (Lemma~\ref{lem:self}). The main result of the paper is Theorem~\ref{thm:main}: it generalizes Theorem~\ref{prop:main1} to the case of piecewise smooth periodic surfaces. Concluding remarks follow last.

\section{Heuristic theory for smooth graphs}
\subsection{Definitions}
\begin{definition}
Let $\Omega$ be a non-empty connected open set. We call \emph{surface} a piecewise smooth mapping
\[
    \begin{split}
        \ten x: \R^2\supset \bar{\Omega}&\to\R^3\\
        (u,v) &\mapsto \ten x(u,v)
    \end{split}
\]
such that the partial derivatives $\ten x_1\equiv\partial\ten x/\partial u$ and $\ten x_2\equiv\partial\ten x/\partial v$ are linearly independent wherever they exist.
\end{definition}

This, certainly reductive, definition of a surface shall be enough for the present purposes. A surface can have discontinuities in its tangent plane that correspond to creases. Note however that continuity is implied by piecewise smoothness. We reserve the term \emph{piecewise differentiable} for mappings that are not necessarily continuous.

\begin{definition}
Let $\ten x$ be a surface and $\dot{\ten x}$ be a piecewise smooth mapping that is smooth wherever $\ten x$ is smooth. The \emph{infinitesimal strain} $\gt\str$ of $\ten x$ deflected by $\dot{\ten x}$ is the $2\times 2$ matrix of components
\[  
    \str_{\mu\nu} \equiv \frac{1}{2}\left(\scalar{\dot{\ten x}_\mu,\ten x_\nu}+\scalar{\dot{\ten x}_\nu,\ten x_\mu}\right).
\]
Then, $\dot{\ten x}$ is an \emph{infinitesimal isometry} if $\gt\str$ vanishes everywhere.
\end{definition}

Thus, admissible deformations $\dot{\ten x}$ are restricted to ones that have the same regularity as $\ten x$: they cannot create new creases.

\begin{definition}
To each infinitesimal isometry $\dot{\ten x}$ of a surface $\ten x$, there corresponds a unique~\cite{Spivak1999a} field $\ten w$ such that
\[\label{eq:infIso}
    \dot{\ten x}_\mu = \ten w\wedge\ten x_\mu.
\]
We call $\ten w$ the field of \emph{infinitesimal rotations} of $\dot{\ten x}$. If $\ten w$ is constant, then $\dot{\ten x}$ is the initial velocity of a euclidean motion and is \emph{trivial}; otherwise, it is \emph{non-trivial}.
\end{definition}

We are mainly concerned with periodic surfaces.

\begin{definition}\label{def:per}
Let $R\subset \R^2$ be a non-empty rectangle. A surface $\ten x$ is $R$-\emph{periodic} if it takes the form
\[
\ten x(u,v) = u\ten p_1 + v\ten p_2 + \tilde{\ten x}(u,v)
\]
where $\tilde{\ten x}$ is a piecewise smooth periodic mapping defined over $\R^2$ with unit cell $R$, and $\ten p_1$ and $\ten p_2$ are two linearly independent vectors in $\R^3$. Then, the $\ten p_\mu$ are the \emph{support vectors} of $\ten x$, mapping $\tilde{\ten x}$ is the \emph{corrugation} and
\[
    I_{\mu\nu} = \scalar{\ten p_\mu,\ten p_\nu}
\]
defines the \emph{effective metric} $\ten I$ of $\ten x$.
\end{definition}

Hereafter, periodicity is always to be understood relative to a particular unit cell; thus, $R$-periodic will be shortened to periodic. Furthermore, in this first part, we focus on smooth periodic graphs. Cases where the graph is piecewise smooth, rather than smooth, are better dealt with in the context of the general theory of the next part.

\begin{definition}
The \emph{graph} $\ten x$ of a piecewise smooth function $z$ is the surface given by
\[
    \ten x(u,v) = (u,v,z(u,v)).
\]
The function $z$ is the \emph{profile} of $\ten x$. The graph is \emph{periodic} (resp. \emph{smooth}) if the profile is periodic (resp. smooth).
\end{definition}

\subsection{Motivation}
Let $\ten x$ be a smooth periodic graph of profile~$z$ and let
\[
    \hat{\ten x} = (\hat x,\hat y,\hat z)
\]
be an infinitesimal isometry of $\ten x$. Then, it is well-known that $\hat z$ is solution to the linear PDE
\[\label{eq:MA}
    \dif_z\hat z \equiv z_{11}\hat z_{22} + z_{22}\hat z_{11} - 2z_{12}\hat z_{12} = 0,
\]
where indices denote partial derivatives. The above equation translates the fact that Gaussian curvature is an isometric invariant. The differential operator $\dif_z$ is the Monge-Ampère operator linearized at $z$; see~\cite{Audoly2010}. Equation~\eqref{eq:MA} admits several trivial solutions. In particular, a constant $\hat z$ corresponds to a vertical translation and a linear $\hat z$ corresponds to an infinitesimal rotation. Here, we are interested in solutions of the form
\[
    \hat z(u,v) = \tilde z(u,v) + q(u,v) \equiv \tilde z(u,v) +  \frac{1}{2}e u^2 + f uv + \frac{1}{2}g v^2,
\]
where $\tilde z$ is periodic and $q$ is a quadratic form of coefficients $(e,f,g)$.
The quadratic form $q$ describes an effective bending of surface $\ten x$ or, in other words, how $\ten x$ bends and twists ``on average''. Specifically, let
\[
\ten x^\epsilon(u,v) = \epsilon\ten x(u/\epsilon,v/\epsilon)
\]
define a surface with refined corrugations for $0<\epsilon<1$. Then, an infinitesimal isometry of $\ten x^\epsilon$ is given by
\[
     \hat z^\epsilon(u,v) \equiv \epsilon^2\hat z(u/\epsilon,v/\epsilon) = \epsilon^2\tilde z(u/\epsilon,v/\epsilon) + \epsilon^2 q(u/\epsilon,v/\epsilon) = \epsilon^2\tilde z(u/\epsilon,v/\epsilon) + q(u,v),
\]
since $q$ is quadratic. Thus, should $\hat z$ exist, then $q$ is an infinitesimal isometry of $\ten x^\epsilon$ up to an error of order $\epsilon^2$; this is what we will refer to as an effective infinitesimal isometry. Indeed, $\norm{\hat z^\epsilon-q}_\infty = O(\epsilon^2)$ because $\tilde z$ is continuous periodic and therefore bounded.

\begin{definition}
    A quadratic form $q$ is an \emph{effective} infinitesimal isometry of a smooth periodic graph of profile~$z$ if there exists a smooth periodic mapping $\tilde z$ such that $\dif_z(\tilde z+q)=0$. In that case, $\tilde z$ is the corrector of $q$.
\end{definition}

Now it is desirable to characterize the effective infinitesimal isometries $q$ without referring to the corrector $\tilde z$. Hereafter, a necessary condition on $q$ is obtained. It relies on a self-adjointness property of the linearized Monge-Ampère operator.

\subsection{Self-adjointness}
The equation to be solved, namely $\dif_z(\tilde z+q)=0$, amounts to solving the linearized Monge-Ampère equation with a non-zero right-hand side. Typically, for a solution to exist, the right hand side must satisfy a compatibility condition as in the Fredholm alternative. Here, the compatibility condition says that the right hand side, function of $q$, must be orthogonal to all periodic solutions~$\dot z$ of $\dif_z\dot z=0$.

\begin{definition}
    An infinitesimal isometry $\dot{\ten x}$ of a smooth periodic graph $\ten x$ is \emph{periodic} if it is of the form $\dot{\ten x} = (\dot x, \dot y, \dot z)$
    where $\dot z$ is periodic, and $\dot x$ and $\dot y$ are periodic modulo linear forms. Then, the \emph{effective} infinitesimal strain $\ten E$ is the $2\times 2$ matrix of coefficients
    \[
        E_{\mu\nu} = \frac{\scalar{\int \dot{\ten x}_\mu,\int \ten x_\nu} + \scalar{\int \dot{\ten x}_\nu,\int \ten x_\mu}}{2}.
    \]
\end{definition}

Above and throughout, the non-annotated symbol of integration $\int$ refers to the \emph{mean value} over the unit cell $R$:
\[
    \int  \equiv \frac{1}{\text{Area}(R)}\int_R \dd\xi_1\dd\xi_2.
\]

The effective infinitesimal strain $\ten E$ quantifies how the unit cell, as immersed in $\R^3$, deforms. In other words, it defines the effective stretch and shear of and among the directions of periodicity, as immersed in $\R^3$. In particular, in reference to the support vectors $\ten p_\mu$ (Definition~\ref{def:per}) and to their deflections $\dot{\ten p}_\mu$, the effective infinitesimal strain is
\[
    E_{\mu\nu} = \frac{\scalar{\dot{\ten p}_\mu,\ten p_\nu}+\scalar{\dot{\ten p}_\nu,\ten p_\mu}}{2}.
\]
Alternatively, the following form of $\ten E$ will prove particularly useful. It involves the gradient operator $\gt\nabla$.

\begin{lemma}\label{lem:E}
Let $\dot{\ten x}=(\dot x,\dot y, \dot z)$ be a periodic infinitesimal isometry of a smooth periodic graph $\ten x$ of profile~$z$. Then, the effective infinitesimal strain is
\[
    \ten E = -\frac{1}{2}\int \left(\gt\nabla \dot z\gt\nabla^T z + \gt\nabla z\gt\nabla^T \dot z\right).
\]
\end{lemma}
\begin{proof}
Set the infinitesimal strain to $0$ to find
\[
\dot x_1 + z_1\dot z_1 = 0,\quad
\dot y_2 + z_2\dot z_2 = 0,\quad
\dot x_2 + \dot y_1  + z_1\dot z_2 + z_2\dot z_1 = 0.
\]
Note that the derivative of a smooth periodic function has zero average. In particular,
\[
\int z_\mu = \int \dot z_\mu = 0.
\]
Therefore,
\[
\begin{split}
    E_{11} &= \scalar{\int \dot{\ten x}_1,\int\ten x_1} = \int \dot x_1 = -\int z_1\dot z_1,\\
    E_{22} &= \scalar{\int \dot{\ten x}_2,\int\ten x_2} = \int \dot y_2 = -\int z_2\dot z_2,\\
    2E_{12} &= \scalar{\int \dot{\ten x}_1,\int\ten x_2}+\scalar{\int \dot{\ten x}_2,\int\ten x_1} = \int \dot x_2 + \int \dot y_1 = -\int z_1\dot z_2-\int z_2\dot z_1,
\end{split}
\]
as claimed.
\end{proof}

It is time to state and prove the anticipated self-adjointness property of the linearized Monge-Ampère operator.
\begin{lemma}\label{lem:self1}
    Let $z$, $\dot z$ and $\tilde z$ be real smooth periodic functions defined over $\R^2$. Then,
    \[
    \int \hat z\dif_z \dot z = \int \dot z\dif_z \hat z.
    \]
\end{lemma}
\begin{proof}
    Use Schwarz theorem to rewrite
    \[
    \dif_z \dot z = \div\left(\adj(\ten H_z)\gt\nabla \dot z\right)
    \]
    where $\div$ is the divergence operator and $\adj(\ten H_f)$ is the adjugate matrix of the Hessian matrix of~$z$. Then apply the divergence theorem twice in conjunction with periodicity and the symmetry of $\adj(\ten H_z)$.
\end{proof}

\subsection{Main result for smooth graphs}
\begin{theorem}\label{prop:main1}
    Let $\ten x$ be a smooth periodic graph that admits a periodic infinitesimal isometry of effective infinitesimal strain $\ten E$. Let $q$ be an effective infinitesimal isometry of $\ten x$ of coefficient $(e,f,g)$. Then
    \[
        \adj(\ten H_q) \cdot \ten E = 0.
    \]
    Explicitly,
    \[
     E_{11}g + E_{22}e - 2E_{12}f = 0.
    \]
\end{theorem}

\begin{proof}
    Let $z$ be the profile of $\ten x$, let $\tilde z$ be the corrector of $\ten q$, and let $\dot z$ be the third component of the periodic infinitesimal isometry corresponding to $\ten E$. Then,
    \[
    \begin{aligned}
      0 &= \int \dot z \dif_z(\tilde z + q) &&\text{since} \quad \dif_z(\tilde z+q)=0\\
        &= \int \tilde z \dif_z \dot z + \int \dot z\dif_zq &&\text{by self-adjointness}\\
        &= \int \dot z\dif_zq &&\text{since} \quad \dif_z\dot z = 0\\
        &= \int \dot z\dif_q z &&\text{since}\quad \dif_a b = \dif_b a\\
        &= -\int \scalar{\gt\nabla\dot z,\adj(\ten H_q)\gt\nabla z} &&\text{by the divergence theorem} \\
        &= -\adj(\ten H_q)\cdot\int \gt\nabla\dot z\gt\nabla^T z &&\text{since $\ten H_q$ is constant}\\
        &= \adj(\ten H_q)\cdot \ten E &&\text{by symmetry and Lemma~\ref{lem:E}}
    \end{aligned}
    \]
\end{proof}

The above theorem clarifies what quadratic surfaces can be embraced, on average, upon an infinitesimal isometric deformation of a periodic smooth graph. Remarkably, the more flexible the surface is in the plane, in the sense that it admits multiple effective infinitesimal strains, the more rigid it is out-of-plane, in the sense that the embraced surfaces have curvatures and torsion that must satisfy multiple constraints of the kind imposed by the theorem. That being said, ``silent'' infinitesimal periodic isometries, i.e., ones with $\ten E=\ten 0$, are inconsequential. To interpret the constraint, it is best to choose the coordinate axes so as to coincide with the axes of principal effective infinitesimal strains. In that case, the constraint simplifies into
\[
E'_{11}g' + E'_{22}e' = 0, \quad \text{or} \quad
-\frac{E'_{11}}{E'_{22}} = \frac{e'}{g'},
\]
if algebra permits. The left hand side is a ratio of relative effective extensions and defines an effective Poisson's ratio; the right hand side is a ratio of effective normal curvatures in the directions of principal effective infinitesimal strains. The theorem proves that the two ratios are equal. Qualitatively, if a smooth periodic graph shrinks laterally as it stretches longitudinally, all inextensibly, then, should it bend on average, it bends into a dome; if, by contrast, it extends laterally then it bends into a saddle. Put succinctly: \emph{auxetic graphs bend anticlastically; anauxetic graphs bend synclastically}.

Another interpretation, in the style of classical geometers, goes as follows: the effective infinitesimal strain defines a paraboloid $p: (u,v) \mapsto E_{11} u^2/2+E_{12} uv+E_{22}v^2/2$ that is either elliptic (case $\det \ten E > 0$), cylindrical (case $\det \ten E = 0$) or hyperbolic (case $\det \ten E < 0$). Similarly, $q$ defines a paraboloid whose nature depends on the sign of $\det \ten H_q$. With that in mind, the constraint reads: $\dif_{p}q=0$, or equivalently $\dif_q p=0$. That is: the paraboloid $p$ admits $q$ as an infinitesimal bending, or equivalently paraboloid $q$ admits $p$ as an infinitesimal bending. This, again, reduces to the invariance of Gaussian curvature of the paraboloids $p$ and $q$ this time. In fact, the Gaussian curvature of paraboloid $q$ is $\det \ten H_q$ and its invariance under deflection $p$ reads
\[
    \det(\ten H_q + t\ten E) = \det \ten H_q + O(t^2),
\]
which, using Jacobi's formula, is the proven constraint. Now say that the paraboloid $p$ is elliptic; then, to bend it isometrically, one of its principal curvatures must increase and the other decrease, meaning that deflection $q$ describes a saddle. If paraboloid $p$ is hyperbolic, both \emph{signed} principal curvatures increase or both decrease, since now the Gaussian curvature of $p$ is negative; thus, deflection $q$ describes a dome. Last, if paraboloid $p$ is cylindrical, one of its principal curvatures stays null, meaning that deflection $q$ describes a cylinder as well.

\subsection{Example: graphs of translation}
Let $a$ and $b$ be two smooth periodic functions. Let $\ten x$ be the periodic graph defined by
\[
    \ten x(u,v) = (u,v,a(u)+b(v)).
\]
This is a \emph{graph of translation}, i.e., a surface swept by the translation of one curve, the \emph{profile}, along another, the \emph{path}. Surface $\ten x$ admits a periodic infinitesimal isometric deformation~\cite{Bianchi1878, izmestiev2023, Nassar2023}, namely
\[
\dot{\ten x}(u,v) = \left(-\int^u a'^2,\int^v b'^2,a(u)-b(v)\right).
\]
The effective infinitesimal strain is
\[
\ten E = \begin{bmatrix}
-\int a'^2 & 0 \\ 0 & \int b'^2
\end{bmatrix}.
\]
Let $(e,f,g)$ be the coefficients of an effective infinitesimal isometry of $\ten x$. Then, by Theorem~\ref{prop:main1},
\[
    \frac{e}{g} = \frac{\int a'^2}{\int b'^2}.
\]
This proves the following, loosely stated, corollary.

\begin{corollary}\label{cor:trans1}
    Smooth periodic graphs of translation effectively bend synclastically.
\end{corollary}

\begin{example}\label{exa:ex}
    Suppose $a=b$, then $e/g=1$ meaning that periodic smooth graphs of translation with identical path and profile bend, on average, ``equi-synclastically''. To verify the pertinence of this proposition for thin elastic shells, a series of finite element simulations were carried out. The model is a shell whose midsurface is given by $a=b=\cos$ and is defined over a range containing an array of $5\times 5$ unit cells. On two edges of the boundary, $\{u=0\}$ and $\{v=0\}$, symmetry boundary conditions are imposed; one edge, $\{u=10\pi\}$, is left free; on the last edge, $\{v=10\pi\}$, a uniform moment about axis $(1,0,0)$ is imposed. The deflections are probed at the locations $\{((2j+1)\pi,0,-1)\}_{j=0\dots 4}$ and $\{(0,(2j+1)\pi, -1)\}_{j=0\dots 4}$ and are fitted with two parabolas whose curvatures are interpreted as $e$ and $g$. An error $\delta = |e/g-1|$ is computed for decreasing thicknesses $h$; the results are reported in Table~\ref{tab:error} and show that the error decreases at least like $h$. Moreover, the numerical results are robust against changes in the material properties. This example shows that the purely geometric theory proposed here of effective isometric deformations is relevant to the deformation of thin elastic shells. Simulations used the Ansys software and employed a structured mesh with (more than) $10^6$ equal-sized elements of type Shell181; see Supplemental Materials for raw data and further detail.
    \begin{table}[ht!]
        \centering
        \begin{tabular}{c|rrrrr}
             $h$     &  $2.5\times10^{-1}$ & $1.25\times10^{-1}$ & $2.5\times 10^{-2}$ & $2.5\times 10^{-3}$ & $2.5\times 10^{-4}$\\
             $\delta$&  $1.26\times 10^{-1}$ & $3.70\times 10^{-2}$ & $1.66\times10^{-3}$ & $5.48\times 10^{-5}$ & $2.40\times10^{-6}$
        \end{tabular}
        \caption{Error $|e/g-1|$ v.s. thickness $h$. Data reported for dimensions in $\mathrm{cm}$, a Young's modulus of $2.7\times 10^9\mathrm{Pa}$ and a material Poisson's ratio of $0.33$; changing material properties has little influence on the trend of $\delta$.}
        \label{tab:error}
    \end{table}
\end{example}

By contrast, Theorem~\ref{prop:main1} leaves effective torsion $f$ free of constraints. Indeed, any $(0,f,0)$ is an effective infinitesimal isometry of any smooth periodic graph of translation.

\begin{corollary}\label{cor:trans2}
    All smooth periodic graphs of translation admit the quadratic form of coefficients $(0,1,0)$ as an effective infinitesimal isometry.
\end{corollary}
\begin{proof}
    Check that $\dif_z q=0$ for $q(u,v)=uv$ and any $z(u,v)=a(u)+b(v)$ meaning that $q$ is an infinitesimal isometry of any smooth graph of translation. In particular, $q$ is an effective infinitesimal isometry of any periodic smooth graph of translation.
\end{proof}

Thus far, the existence of an effective infinitesimal strain $\ten E$ informed on effective infinitesimal isometries $q$. But Theorem~\ref{prop:main1} goes both ways.

\begin{corollary}\label{cor:trans3}
    The profile and path of a smooth periodic graph of translation are unshearable.
\end{corollary}
\begin{proof}
    The quadratic form of coefficients $(0,1,0)$ being an effective infinitesimal isometry implies $E_{12}=0$ for any effective infinitesimal strain tensor $\ten E$.
\end{proof}

\section{Asymptotic theory for piecewise smooth surfaces}
The present section extends Theorem~\ref{prop:main1} to effective \emph{finite} isometries and deals with the complications that arise from considering surfaces that are not graphs or are not smooth. 
\subsection{Preliminaries}
A key element of the theory that follows is a generalization of the self-adjointness property of Lemma~\ref{lem:self1} as well as of the operator $\dif_z$ itself. In this generalization, it is more natural to deal with rotations than deflections.

\begin{definition}\label{def:adm}
    A piecewise differentiable field $\ten w$ is an \emph{admissible} field of infinitesimal rotations of a surface $\ten x$ if $s\mapsto \ten w\wedge\dd\ten x/\dd s$ is single-valued for any $s\mapsto (u(s),v(s))$ that parametrizes a line of discontinuity in the tangent plane of~$\ten x$.
\end{definition}

Simply put, admissibility defines candidate infinitesimal rotations $\ten w$ that could produce infinitesimal isometries. In particular, candidate rotations must maintain continuity across crease lines. This implies that rotations, on either side of the crease, must rotate the crease in the same fashion. In other words, $\ten w\wedge\dd\ten x/\dd s$ must be single-valued.

\begin{lemma}
    If $\ten w$ is the field of infinitesimal rotation of an infinitesimal isometry $\dot{\ten x}$ of a surface $\ten x$, then
\[
    \dar_{\ten x}\ten w \equiv \ten w_v\wedge\ten x_u - \ten w_u\wedge\ten x_v = \ten 0,
\]
and $\ten w$ is admissible.
\end{lemma}
\begin{proof}
    Equation~\eqref{eq:infIso} implies $\ten w$ is piecewise differentiable because $\ten x$ and $\dot{\ten x}$ are. Furthermore, the continuity of $\dot{\ten x}$ implies the admissibility of $\ten w$ by the, here admitted, Hadamard jump condition. Last, the same equation yields
    \[
    \ten 0 = \dot{\ten x}_{uv}-\dot{\ten x}_{vu} = \dar_{\ten x}\ten w.
    \]
\end{proof}

The $\dif$ was for Monge. Here, the $\dar$ is for Darboux. Like $\dif$, operator $\dar$ is self-adjoint.

\begin{lemma}\label{lem:self}
    Let $\ten x$ be a periodic surface. Then,
    \[
        \int\scalar{\gt\omega,\dar_{\ten x}\ten w} = \int\scalar{\ten w,\dar_{\ten x}\gt\omega},
    \]
    for any $\gt\omega$ and $\ten w$ that are periodic and admissible
\end{lemma}

\begin{proof}
    Let $\{R_i\}_{1\leq i\leq n}$ be a finite set of disjoint non-empty open connected sets such that $\ten x$ is smooth over $\bar{R_i}$ and such that $\cup_i\bar{R}_i=\bar R$, where $R$ is the unit cell. Let $\partial R_{ij}=\bar{R}_i\cap(\bar{R}_j+\mathcal{R})$, where $\mathcal{R}$ is the periodicity lattice. Let $s\mapsto(u(s),v(s))$ parametrize one of these intersections and let the brackets $[\cdot]$ denote the jump in any quantity across the intersection. Then,
    \[\label{eq:jump}
        \begin{aligned}\relax
            [\scalar{\gt\omega,\ten w\wedge\dd\ten x/\dd s}]
            &= \scalar{[\gt\omega],\ten w\wedge\dd\ten x/\dd s} &&\text{since $\ten w$ is admissible}\\
            &= \scalar{\ten w,\dd\ten x/\dd s\wedge [\gt\omega]} &&\text{by permutation symmetry}\\
            &= \scalar{\ten w,[\dd\ten x/\dd s\wedge \gt\omega]} &&\text{by continuity of $\ten x$}\\
            &= 0 &&\text{since $\gt\omega$ is admissible}.
        \end{aligned}
    \]
    Now write
    \[
        \begin{aligned}
            &\int_R\scalar{\gt\omega,\dar_{\ten x}\ten w} \\
            =& \int_R\scalar{\gt\omega,\ten w_v\wedge\ten x_u - \ten w_u\wedge\ten x_v} &&\text{by definition}\\
            =& \int_R\scalar{\gt\omega,(\ten w\wedge\ten x_u)_v - (\ten w\wedge\ten x_v)_u} &&\text{by Schwarz theorem}\\
            =& \sum_{i}\oint_{\partial R_{i}}\scalar{\gt\omega,\ten w\wedge\dd\ten x/\dd s} - \int_R\left(\scalar{\gt\omega_v,\ten w\wedge\ten x_u}-\scalar{\gt\omega_u, \ten w\wedge\ten x_v}\right) &&\text{by the divergence theorem}\\
            =& \sum_{i}\oint_{\partial R_{i}}\scalar{\gt\omega,\ten w\wedge\dd\ten x/\dd s} + \int_R\scalar{\ten w,\dar_{\ten x}\gt\omega} &&\text{by permutation symmetry}\\
            =& \sum_{i<j}\int_{\partial R_{ij}}[\scalar{\gt\omega,\ten w\wedge\dd\ten x/\dd s}] + \int_R\scalar{\ten w,\dar_{\ten x}\gt\omega} &&\text{since $\partial R_i=\cup_j\partial R_{ij}$}\\
            =& \int_R\scalar{\ten w,\dar_{\ten x}\gt\omega} &&\text{by equation~\eqref{eq:jump}}.
        \end{aligned}
    \]
\end{proof}

Last, it is worthwhile to generalize the notion of periodicity.

\begin{definition}
    An infinitesimal isometry $\dot{\ten x}$ of a periodic surface $\ten x$ is \emph{periodic} if it is of the form
    \[
        \dot{\ten x}(u,v) = u\dot{\ten p}_1 + v\dot{\ten p}_2 + \dot{\tilde{\ten x}}(u,v)
    \]
    where $\dot{\tilde{\ten x}}$ is periodic, and $\dot{\ten p}_1$ and $\dot{\ten p}_2$ are two vectors in $\R^3$. Then, the \emph{effective} infinitesimal strain $\ten E$ is the $2\times 2$ matrix of coefficients
    \[
        E_{\mu\nu} = \frac{\scalar{\int \dot{\ten x}_\mu,\int \ten x_\nu} + \scalar{\int \dot{\ten x}_\nu,\int \ten x_\mu}}{2}.
    \]
\end{definition}

\subsection{Asymptotic expansions}
The idea as before is that an effective (finite) isometry preserves lengths to within an error in~$O(\epsilon^2)$ where $\epsilon$ is the size of the corrugations. The following definitions codify this notion.

\begin{definition}
    Two surfaces $\ten x$ and $\ten y$ are \emph{isometric} if they have the same metric, namely
    \[
        \scalar{\ten x_\mu,\ten x_\nu} = \scalar{\ten y_\mu,\ten y_\nu},
    \]
    and if $\ten y$ is smooth wherever~$\ten x$ is. If $\ten y = \ten R\ten x + \ten T$ for some constant rotation matrix~$\ten R$ and constant vector~$\ten T$, then the isometry is \emph{trivial}; otherwise, it is \emph{non-trivial}.
\end{definition}

\begin{definition}
    Let $\ten x$ be a periodic surface and let, for $\epsilon>0$,
    \[
    \begin{split}
        \ten x^\epsilon:\bar{\Omega}&\to\R^3\\
        (u,v)&\mapsto \epsilon\ten x(u/\epsilon,v/\epsilon).
    \end{split}
    \]
    Let $\ten y^\epsilon$ be a surface that is smooth wherever $\ten x^\epsilon$ is and such that
    \[\label{eq:effbend}
        \begin{split}
        \scalar{\ten x^\epsilon_\mu,\ten x^\epsilon_\nu} - \scalar{\ten y^\epsilon_\mu,\ten y^\epsilon_\nu} & = O(\epsilon^2),\\
        \lim_{\epsilon\to 0}\ten y^\epsilon &= \ten Y,
        \end{split}
    \]
    where $\ten Y:\bar{\Omega}\to\R^3$ is a surface. Then, $\ten Y$ is an \emph{effective isometry} of $\ten x$ over $\bar{\Omega}$.
\end{definition}

Here, convergence is understood in the sense of $\lVert\cdot\rVert
_\infty$ over $\bar\Omega$. Now in the spirit of the theory of homogenization, it is desirable to characterize the effective isometries $\ten Y$ without referring to the sequences $\ten y^\epsilon$. This is pursued hereafter for a specific class of effective isometries that are the limit of an ansatz of the form
\[\label{eq:ansatz}
    \begin{split}
        \ten y^\epsilon(u,v) &= \ten z^\epsilon(u,v,u/\epsilon,v/\epsilon),\\
        \ten z^\epsilon(U,V,u,v) &= \ten Y(U,V) + \epsilon \ten z^1(U,V,u,v) + \epsilon^2 \ten z^2(U,V,u,v),
    \end{split}
\]
where $\ten z^\epsilon$, $\ten z^1$ and $\ten z^2$ are piecewise smooth periodic mappings in $(u,v)$ and piecewise smooth in~$(U,V)$. The ansatz is rather intuitive and typical of problems where a PDE involves periodic rapidly oscillating coefficients, namely the $\ten x^\epsilon_\mu$. That being said, we do not claim that all effective isometries are limit of an ansatz of this form.

Now substitute the ansatz~\eqref{eq:ansatz} into~\eqref{eq:effbend} to obtain
    \begin{multline}\label{eq:expansion}
     \scalar{\ten x_\mu,\ten x_\nu} - \scalar{\ten Y_M + \ten z^1_\mu, \ten Y_N + \ten z^1_\nu}
     \\- \epsilon\left(\scalar{\ten Y_M + \ten z^1_\mu, \ten z^1_N+\ten z^2_\nu} + \scalar{\ten z^1_M+\ten z^2_\mu, \ten Y_N + \ten z^1_\nu}\right) = O(\epsilon^2).
     \end{multline}
Therein, and from now on, it is understood that $M=\mu$ and $N=\nu$ albeit $(M,N)$ (resp., $(\mu,\nu)$) denote derivatives relative to $(U,V)$ (resp., $(u,v)$). Thanks to the following lemma, it is legitimate to ``balance'' terms of the same order in $\epsilon$ at any point $(U,V,u,v)$.

\begin{lemma}\label{lem:TheLemma}
    Let $a(u,u/\epsilon) = O(\epsilon)$ where $a$ is continuous and periodic in its second variable. Then, $a=0$.
\end{lemma}
\begin{proof}
    Admitted.
\end{proof}

This yields, for leading-order terms,
\[\label{eq:order1}
    \scalar{\ten x_\mu,\ten x_\nu} = \scalar{\ten Y_M + \ten z^1_\mu, \ten Y_N + \ten z^1_\nu}
\]
for all $(U,V,u,v)$. Hence, for any $(U,V)$, $\ten b(u,v) \equiv \ten Y_U u + \ten Y_V v + \ten z^1(u,v)$ describes a periodic surface that is isometric to $\ten x$. These isometries could involve trivial ones as well as various non-trivial ones. Moving forward, we focus on a particular case outlined by the following definition.

\begin{definition}\label{def:unimode}
    Let $\ten x$ be a periodic surface that admits a continuously differentiable one-parameter family of isometric periodic surfaces $I\ni t\mapsto \ten b(t)$ where
    \[
        \ten b(t): (u,v)\mapsto u\ten p_1(t) + v\ten p_2(t) + \tilde{\ten b}(t,u,v).
    \]
    Then, a \emph{unimodal} effective isometry $\ten Y$ of \emph{mode} $\ten b$ is the limit of an ansatz of the form
        \begin{align}
        \ten y^\epsilon(u,v) &= \ten z^\epsilon(u,v,u/\epsilon,v/\epsilon),\nonumber\\
        \ten z^\epsilon(U,V,u,v) &= \ten Y(U,V) + \epsilon \ten z^1(U,V,u,v) + \epsilon^2 \ten z^2(U,V,u,v),\\
        \ten z^1(U,V,u,v) &= \ten R(U,V)\ten b(\theta(U,V),u,v)-\ten Y_U(U,V)u -\ten Y_V(U,V)v + \ten T(U,V),\nonumber
    \end{align}
    for some piecewise differentiable functions
    \[
        \theta: \bar\Omega \to I,\qquad
        \ten R: \bar\Omega \to \mathcal{O}_3(\R),\qquad
        \ten T: \bar\Omega \to \R^3.
    \]
\end{definition}

Thus, a unimodal effective isometry $\ten Y$ of a periodic surface $\ten x$ is the limit of deformations $\ten y^\epsilon$ that ``modulate'' in space one given mode of isometric deformation of $\ten x$, namely $\ten b$. Other effective isometries could involve multiple modes or refute the proposed ansatz altogether.

\subsection{Main result for piecewise smooth surfaces}
\begin{theorem}\label{thm:main}
    Let $\ten Y$ be a unimodal effective isometry of mode
    \[
        \ten b(t): (u,v)\mapsto u\ten p_1(t) + v\ten p_2(t) + \tilde{\ten b}(t,u,v),
    \]
    of a periodic surface $\ten x$. Let $\ten I$ be the effective metric of $\ten b$ and let $\ten E\equiv\dot{\ten I}/2$ be the effective infinitesimal strain of $\dot{\ten b}\equiv\dd \ten b/\dd t$. Then, there exists a piecewise differentiable function $\theta: \bar\Omega \to I$ such that
    \[\label{eq:Y1}
        \scalar{\ten Y_M,\ten Y_N} = I_{\mu\nu}(\theta)
    \]
    and
    \[\label{eq:Y2}
        E_{22}(\theta)\ten Y_{11} + E_{11}(\theta)\ten Y_{22}
        -2E_{12}(\theta)\ten Y_{12}
        = \ten 0.
    \]
\end{theorem}

Before diving into the proof, a few comments are helpful. Since $\ten b$ is a continuously differentiable family of isometries, the velocity field $\dd \ten b/\dd t$ is, at each $t$, an infinitesimal isometry of $\ten b(t)$. Then, the content of the theorem is that the finite isometry provides the metric of the effective isometry $\ten Y$, i.e., the effective metric, whereas the infinitesimal isometry provides, by the same logic of Theorem~\ref{prop:main1}, a constraint on the effective curvatures. The first part of the theorem falls directly from the definition of $\ten Y$; the challenge resides in proving the second part. This involves, as in the case of smooth graphs, a self-adjointness property that is necessary for the existence of a solution to the equation of infinitesimal isometries with a non-zero right hand side. The self-adjointness property is the one stated in Lemma~\ref{lem:self}.

\begin{proof}
    Equation~\eqref{eq:Y1} follows from
    \[
        \begin{aligned}
            \ten Y_M
            &= \int\ten Y_M && \text{since $\ten Y$ is $(u,v)$-independent}\\
            &= \int\ten R\ten b_\mu - \int\ten z^1_\mu && \text{by Definition~\ref{def:unimode}}\\
            &= \int\ten R\ten b_\mu && \text{by periodicity of $\ten z^1$}\\
            &= \ten R\int\ten b_\mu && \text{since $\ten R$ is $(u,v)$-independent}\\
            &= \ten R\ten p_\mu && \text{by periodicity of $\tilde{\ten b}$},
        \end{aligned}
    \]
    from Definition~\ref{def:per}, and from the fact that $\ten R$ is a rotation.
    
    Back to~\eqref{eq:expansion}, balancing the terms linear in $\epsilon$ yields
    \[
        \scalar{\ten R\ten b_\mu, \ten z^1_N+\ten z^2_\nu} + \scalar{\ten z^1_M+\ten z^2_\mu, \ten R \ten b_\nu} = 0.
    \]
    Thus, there exists a vector field $\ten w$ such that
    \[\label{eq:w}
        \ten z^1_M+\ten z^2_\mu = \ten w \wedge \ten R\ten b_\mu.
    \]
    Note that, relative to the surface $\ten R\ten b$, field $\ten w$ is admissible in the sense of Definition~\ref{def:adm} since $\ten z^1_M$ and $\ten z^2$ are continuous. Furthermore, $\ten z^2_{\mu\nu}=\ten z^2_{\nu\mu}$ entails
    \[
        \ten z^1_{Uv} - \ten z^1_{Vu} = (\ten w \wedge \ten R\ten b_u)_v - (\ten w \wedge \ten R\ten b_v)_u = \dar_{\ten R\ten b}\ten w.
    \]
    Similarly, let $\gt\omega$ be the infinitesimal rotation vector of the infinitesimal isometry $\dot{\ten b}$ of surface $\ten b$; then $\ten R\gt\omega$ is the infinitesimal rotation vector of the infinitesimal isometry $\ten R\dot{\ten b}$ of surface $\ten R\ten b$. In particular, $\ten R\dot{\gt\omega}$ is admissible and
    \[
        \dar_{\ten R\ten b}\ten R\gt\omega = \ten 0.
    \]
    Then, Lemma~\ref{lem:self} yields
    \[
        0 = \int\scalar{\ten w,\dar_{\ten R\ten b}\ten R\gt\omega} = \int\scalar{\ten R\gt\omega,\dar_{\ten R\ten b}\ten w} = \int{\scalar{\ten R\gt\omega, \ten z^1_{Uv} - \ten z^1_{Vu}}}.
    \]
    Upon substituting in the expressions of $\ten z^1_u$ and $\ten z^1_v$, it comes that
    \[
        \int{\scalar{\ten R\gt\omega, \ten R_U\ten b_v + \theta_U\ten R\dot{\ten b}_v - 
        \ten R_V\ten b_u - \theta_V\ten R\dot{\ten b}_u}} = 0
    \]
    and that
    \[
        \int{\scalar{\ten R\gt\omega, \ten R_U\ten b_v - 
        \ten R_V\ten b_u}} = 0
    \]
    because $\dot{\ten b}_\mu=\gt\omega\wedge\ten b_\mu\perp\gt\omega$. Given that $\ten R^T\ten R_U$ and $\ten R^T\ten R_V$ are skew-symmetric, there exist vectors $\ten\Omega^M$ such that
    \[
    \int{\scalar{\gt\omega, \ten\Omega^1\wedge\ten b_v - 
        \ten\Omega^2\wedge\ten b_u}} = 0.
    \]
    Re-arrange into
    \[
    \scalar{\ten\Omega^1,\int\gt\omega\wedge\ten b_v} - 
        \scalar{\ten\Omega^2,\int\gt\omega\wedge\ten b_u} = 0
    \]
    to find that
    \[\label{eq:comp3}
        \scalar{\ten\Omega^1,\dot{\ten p}_2}=\scalar{\ten\Omega^2,\dot{\ten p}_1}.
    \]
    The last part of the proof aims to recast the above equation into one that involves $\ten Y$ directly rather than the $\ten\Omega$'s.

    Derive $\ten Y_M = \ten R\ten p_\mu$ to get
    \[
        \ten Y_{MN} = \ten R_N\ten p_\mu + \theta_N\ten R\dot{\ten p}_\mu
    \]
    and multiply by $\ten R^T$ to get
    \[
        \ten R^T\ten Y_{MN} = \ten\Omega^N\wedge\ten p_\mu + \theta_N\dot{\ten p}_\mu,
    \]
    and, in particular,
    \[\label{eq:YUV}
        \ten R^T\ten Y_{UV} = \ten\Omega^2\wedge\ten p_1 + \theta_V\dot{\ten p}_1 = \ten\Omega^1\wedge\ten p_2 + \theta_U\dot{\ten p}_2.
    \]
    Thus,
    \[
    \begin{split}
       \scalar{\ten R^T\ten Y_{UU},\ten p_1} &= \theta_U\scalar{\dot{\ten p}_1,\ten p_1},\\
       \scalar{\ten R^T\ten Y_{UV},\ten p_1} &= \theta_V\scalar{\dot{\ten p}_1,\ten p_1},\\
       \scalar{\ten R^T\ten Y_{VV},\ten p_1} &= \scalar{\ten\Omega^2\wedge\ten p_2,\ten p_1}+\theta_V\scalar{\dot{\ten p}_2,\ten p_1}\\
       &= -\scalar{\ten\Omega^2\wedge\ten p_1,\ten p_2}+\theta_V\scalar{\dot{\ten p}_2,\ten p_1}\\
       &= \scalar{\theta_V\dot{\ten p}_1-\theta_U\dot{\ten p}_2,\ten p_2}+\theta_V\scalar{\dot{\ten p}_2,\ten p_1}
    \end{split}
    \]
    where the last equality comes from~\eqref{eq:YUV}. Combine into
    \begin{multline}
    \scalar{\dot{\ten p}_2,\ten p_2}\scalar{\ten R^T\ten Y_{UU},\ten p_1} - (\scalar{\dot{\ten p}_1,\ten p_2}+\scalar{\dot{\ten p}_2,\ten p_1})\scalar{\ten R^T\ten Y_{UV},\ten p_1} + \scalar{\dot{\ten p}_1,\ten p_1}\scalar{\ten R^T\ten Y_{VV},\ten p_1} = 0.
    \end{multline}
    Similarly, it can be shown that
    \begin{multline}
    \scalar{\dot{\ten p}_2,\ten p_2}\scalar{\ten R^T\ten Y_{UU},\ten p_2} - (\scalar{\dot{\ten p}_1,\ten p_2}+\scalar{\dot{\ten p}_2,\ten p_1})\scalar{\ten R^T\ten Y_{UV},\ten p_2} + \scalar{\dot{\ten p}_1,\ten p_1}\scalar{\ten R^T\ten Y_{VV},\ten p_2} = 0.
    \end{multline}
    Finally, consider
    \[\label{eq:NormalProjections}
    \begin{split}
        \scalar{\ten R^T\ten Y_{UU},\dot{\ten p}_1\wedge\dot{\ten p}_2} 
            &= \scalar{\ten\Omega^1\wedge\ten p_1,\dot{\ten p}_1\wedge\dot{\ten p}_2}\\ 
            &= \scalar{\ten\Omega^1,\ten p_1\wedge(\dot{\ten p}_1\wedge\dot{\ten p}_2)}\\
            &= \scalar{\ten\Omega^1,\scalar{\ten p_1,\dot{\ten p}_2}\dot{\ten p}_1-\scalar{\ten p_1,\dot{\ten p}_1}\dot{\ten p}_2}\\
            &= \scalar{\ten p_1,\dot{\ten p}_2}\scalar{\ten\Omega^1,\dot{\ten p}_1}-\scalar{\ten p_1,\dot{\ten p}_1}\scalar{\ten\Omega^1,\dot{\ten p}_2},\\
        \scalar{\ten R^T\ten Y_{VV},\dot{\ten p}_1\wedge\dot{\ten p}_2} 
            &= \scalar{\ten\Omega^2\wedge\ten p_2,\dot{\ten p}_1\wedge\dot{\ten p}_2}\\ 
            &= \scalar{\ten\Omega^2,\ten p_2\wedge(\dot{\ten p}_1\wedge\dot{\ten p}_2)}\\
            &= \scalar{\ten\Omega^2,\scalar{\ten p_2,\dot{\ten p}_2}\dot{\ten p}_1-\scalar{\ten p_2,\dot{\ten p}_1}\dot{\ten p}_2}\\
            &= \scalar{\ten p_2,\dot{\ten p}_2}\scalar{\ten\Omega^2,\dot{\ten p}_1}-\scalar{\ten p_2,\dot{\ten p}_1}\scalar{\ten\Omega^2,\dot{\ten p}_2},\\
        \scalar{\ten R^T\ten Y_{UV},\dot{\ten p}_1\wedge\dot{\ten p}_2} &=
        \scalar{\ten p_1,\dot{\ten p}_2}\scalar{\ten\Omega^2,\dot{\ten p}_1} - \scalar{\ten p_1,\dot{\ten p}_1}\scalar{\ten\Omega^2,\dot{\ten p}_2}\\
        &=
        \scalar{\ten p_2,\dot{\ten p}_2}\scalar{\ten\Omega^1,\dot{\ten p}_1} - \scalar{\ten p_2,\dot{\ten p}_1}\scalar{\ten\Omega^1,\dot{\ten p}_2}.
    \end{split}
    \]
    Then, with~\eqref{eq:comp3} in mind,
    \begin{multline}
    \scalar{\dot{\ten p}_2,\ten p_2}\scalar{\ten R^T\ten Y_{UU},\dot{\ten p}_1\wedge\dot{\ten p}_2} + \scalar{\dot{\ten p}_1,\ten p_1}\scalar{\ten R^T\ten Y_{VV},\dot{\ten p}_1\wedge\dot{\ten p}_2} \\=
    \scalar{\dot{\ten p}_2,\ten p_2}\scalar{\ten p_1,\dot{\ten p}_2}\scalar{\ten\Omega^1,\dot{\ten p}_1}-\scalar{\dot{\ten p}_1,\ten p_1}\scalar{\ten p_2,\dot{\ten p}_1}\scalar{\ten\Omega^2,\dot{\ten p}_2}.
    \end{multline}
    Substituting expressions for $\scalar{\dot{\ten p}_2,\ten p_2}\scalar{\ten\Omega^1,\dot{\ten p}_1}$ and $\scalar{\dot{\ten p}_1,\ten p_1}\scalar{\ten\Omega^2,\dot{\ten p}_2}$ in terms of $\scalar{\ten R^T\ten Y_{UU},\dot{\ten p}_1\wedge\dot{\ten p}_2}$ obtained from~\eqref{eq:NormalProjections}, again with~\eqref{eq:comp3} in mind, leads to
    \begin{multline}
     \scalar{\dot{\ten p}_2,\ten p_2}\scalar{\ten R^T\ten Y_{UU},\dot{\ten p}_1\wedge\dot{\ten p}_2} + \scalar{\dot{\ten p}_1,\ten p_1}\scalar{\ten R^T\ten Y_{VV},\dot{\ten p}_1\wedge\dot{\ten p}_2} \\= (\scalar{\dot{\ten p}_1,\ten p_2}+\scalar{\dot{\ten p}_2,\ten p_1})\scalar{\ten R^T\ten Y_{UV},\dot{\ten p}_1\wedge\dot{\ten p}_2}.
    \end{multline}
Assuming $(\ten p_1,\ten p_2,\dot{\ten p}_1\wedge\dot{\ten p}_2)$ is a basis, it can be concluded that
\[
    \scalar{\dot{\ten p}_2,\ten p_2}\ten Y_{UU} + \scalar{\dot{\ten p}_1,\ten p_1}\ten Y_{VV} = (\scalar{\dot{\ten p}_1,\ten p_2}+\scalar{\dot{\ten p}_2,\ten p_1})\ten Y_{UV}.
\]
Should $\dot{\ten p}_1\wedge\dot{\ten p}_2$ be in the plane of $(\ten p_1,\ten p_2)$, change mode $\ten b$ by superposing a rigid body motion that ensures that $(\ten p_1,\ten p_2,\dot{\ten p}_1\wedge\dot{\ten p}_2)$ is a basis and note that the above equation is invariant under said superposition.
\end{proof}

Theorem~\ref{thm:main} provides two constraints on the unimodal effective isometries of a periodic surface $\ten x$. First, by equation~\eqref{eq:Y1}, the metric of the effective isometry is to be selected among a restricted set of metrics provided by the underlying mode. This is typical of problems where a particular parametrization of a surface is sought, e.g., a conformal parametrization or a Chebyshev net. In such cases, the only constraints that weigh on the second fundamental form are the Gauss-Codazzi-Mainardi equations. The theorem says there is more: the second fundamental form must satisfy an additional algebraic constraint, namely
\[
    \adj(\ten{II})\cdot\ten E(\theta)=0,
\]
where $\ten{II}$ is the second fundamental form of $\ten Y$. This is nothing but equation~\eqref{eq:Y2} projected over the unit normal to $\ten Y$. This constraint affords the same discussion as Theorem~\ref{prop:main1}. In particular: \emph{auxetic surfaces bend anticlastically; anauxetic surfaces bend synclastically}, with the ratio of normal curvatures in the principal directions of effective infinitesimal strain given by the effective Poisson's ratio taken along the same directions. That being said, Theorem~\ref{thm:main} is stronger than Theorem~\ref{prop:main1}: it provides a PDE, namely equation~\eqref{eq:Y2}, that can be solved for $\ten Y$. Remarkably, the PDE is quasilinear and autonomous. Moreover, the type of the PDE is function of the sign of the determinant of the effective infinitesimal strain. In cases where the determinant is positive (resp. negative, null), the PDE is elliptic (resp. hyperbolic, parabolic). The type of the PDE guides appropriate choices of boundary conditions that would guarantee a notion of well-posedness and is therefore critical for applications where the shape of the surface ``in the bulk'' is controlled through boundary input. Morally: to stabilize auxetic surfaces, a closed contour should be pinned; to stabilize anauxetic surfaces, part of the contour should be built-in.

Evidently, the theorem is only useful in cases where a mode $\ten b$ is known. This is thematic in homogenization theory where the use of effective field equations rely on knowledge of solutions to unit cell problems. The theorem is most insightful in cases where mode $\ten b$ is unique (modulo a rigid body motion) as is the case of generic quad-based polyhedral surfaces~\cite{Schief2008}. The theorem also says that the unimodal effective isometries of a trivial mode are to be found among the isometries of the plane: a rather reasonable statement.

The main shortcoming of the theorem is that it provides but a necessary condition on effective isometries. In other words, given a solution $\ten Y$, the theorem does not guarantee the existence of a corresponding sequence $\ten y^\epsilon$. The obvious objection to the converse is that the constraint on $\ten Y$ translates how its second derivatives should be arranged at $(U,V)$ with respect to the infinitesimal isometry $\dot{\ten b}(\theta(U,V))$. Thus, should other infinitesimal isometries $\dot{\ten x}$ be available, there would be other constraints that the second derivatives of $\ten Y$ must satisfy. It is nonetheless tempting to conjecture that Thereom~\ref{thm:main} has a reciprocal in cases where $\dot{\ten b}$ is the unique periodic infinitesimal isometry, modulo rigid body motions.

Many illustrations of the theory presented here exist in the literature. The property that the effective Poisson's ratio is equal to the ratio of normal curvatures in principal directions of strain has been shown to hold, either numerically or analytically, for many quad-based origami tessellations including the Miura-ori~\cite{Schenk2013, Wei2013}, the ``eggbox'' pattern~\cite{Nassar2017a}, the ``morph'' pattern~\cite{Pratapa2019}, and ``zigzag sums''~\cite{Nassar2022}. More recently, this property was numerically shown to hold in a class of curved-crease origami tessellations~\cite{Karami2023}. Particular axisymmetric solutions $\ten Y$ were obtained in~\cite{Nassar2017a,Nassar2017e,Nassar2018b}; see also the numerical scheme of~\cite{Marazzato2023}. Reference~\cite{Nassar2022} verifies, in a particular case, that $\ten Y$ derives from a sequence $y^\epsilon$ that preserve lengths up to an error of order $O(\epsilon^2)$.

\section{Conclusion}
The use of the flexure theory of shells is hindered by the non-existence of non-trivial isometric deformations. As a remedy, and in the context of periodic shells, we propose to extend the available notions of isometric deformations to include \emph{effective isometric deformations} defined as deformations that are inextensional up to first order in the periodicity wavelength. The paper explores this notion and proves two main theorems that provide necessary conditions on the existence of effective isometric deformations for smooth periodic graphs and for general piecewise smooth periodic surfaces, respectively. Morally, the theorems place constraints on how periodic shells can bend on average relative to how they can extend on average. In particular, they show that \emph{auxetic shells bend synclastically} into domes and that \emph{anauxetic shells bend anticlastically} into saddles. Furthermore, the ratio of effective normal curvatures in the principal directions of effective infinitesimal strain is equal to the effective Poisson's ratio taken along these same directions. This property, previously thought to be particular to certain origami tessellations, is here proven to hold for general periodic shells with or without straight and curved creases. 

A cornerstone of the presented theory is a self-adjointness property of the differential operator of infinitesimal isometries restricted to periodic deformations. This property could have other implications on the general theory of isometric deformations that are worthy of investigation. 

Finally, the asymptotic ansatz employed in Theorem~\ref{thm:main} appears to be pertinent when the periodic surface has a non-trivial mode of isometric deformation. For surfaces that only possess trivial modes, other ansatz could be better suited: recent case studies~\cite{Reddy2023} and~\cite{Imada2023} suggest that mesoscales in $\sqrt{\epsilon}$ become relevant.

\section*{Acknowledgments}
This work is funded by NSF CMMI CAREER award no. 2045881. The authors have no conflicts of interest.


\end{document}